\documentclass[graybox]{svmult}

\usepackage{subcaption} 

\usepackage[show]{ed}

\usepackage{type1cm}        
\usepackage{makeidx}         
\usepackage{graphicx}        
\usepackage{multicol}        
\usepackage[bottom]{footmisc}

\usepackage{newtxtext}       %
\usepackage{newtxmath}       

\makeindex             

\usepackage{multirow}
\usepackage{soul}
\usepackage{tabularx,ragged2e}
\newcolumntype{C}{>{\Centering\arraybackslash}X} 
\usepackage{caption}
\usepackage{subcaption}
\usepackage{float}
\usepackage{xcolor}
\usepackage{pdfpages}

\newcommand{\myP}{P}

\newcommand{\beq}{\begin{equation}}
\newcommand{\eeq}{\end{equation}}

\newcommand{\cR}{\mathcal{R}}

\newcommand{\cN}{{\mathcal N}}

\newcommand{\supp}{\mathrm{supp}}

\usepackage{color}
\usepackage{hyperref}
\definecolor{myblue}{rgb}{0,0,0.6}
\definecolor{darkgreen}{rgb}{0,0.5,0}
\hypersetup{colorlinks=true,
linkcolor=myblue,citecolor=myblue,filecolor=myblue,urlcolor=myblue}

\definecolor{escol}{rgb}{0,0.4,0}
\definecolor{sgcol}{rgb}{0,0,0.7}
\definecolor{iggcol}{rgb}{0.5,0,0}
\definecolor{esnewcol}{rgb}{0,0.5,0}

\newcommand{\igg}[1]{{\color{iggcol}{#1}}}

\newcommand{\iggnote}[1]{\ednote{\igg{Ivan  says: #1}}}

\newcommand{\beqs}{\begin{equation*}}
\newcommand{\eeqs}{\end{equation*}}
\newcommand{\bit}{\begin{itemize}}
\newcommand{\eit}{\end{itemize}}
\newcommand{\ben}{\begin{enumerate}}
\newcommand{\een}{\end{enumerate}}
\newcommand{\bal}{\begin{align}}
\newcommand{\eal}{\end{align}}
\newcommand{\bals}{\begin{align*}}
\newcommand{\eals}{\end{align*}}
\newcommand{\bre}{\begin{remark}}
\newcommand{\ere}{\end{remark}}
\newcommand{\bpf}{\begin{proof}}
\newcommand{\epf}{\end{proof}}
\newcommand{\ble}{\begin{lemma}}
\newcommand{\ele}{\end{lemma}}
\newcommand{\bco}{\begin{corollary}}
\newcommand{\eco}{\end{corollary}}
\newcommand{\bex}{\begin{example}}
\newcommand{\eex}{\end{example}}
\newcommand{\bth}{\begin{theorem}}
\newcommand{\enth}{\end{theorem}}

\newcommand{\tfa}{\text{ for all }}

\newcommand{\tand}{\text{ and }}

\newcommand*{\N}[1]{\left\|#1\right\|}

\begin{document}

\title*{Schwarz methods with PMLs for  Helmholtz problems: fast convergence at high frequency}
\titlerunning{Schwarz methods with PMLs for Helmholtz}
\author{J.~Galkowski, S.~Gong, I.G.~Graham, D.~Lafontaine, E.A.~Spence}
\authorrunning{J.~Galkowski, S.~Gong, I.G.~Graham, D.~Lafontaine, E.A.~Spence}

\institute{Jeffrey Galkowski$^1$, Shihua Gong$^2$, Ivan G.~Graham$^3$, David Lafontaine$^4$ and Euan A.~Spence$^3$ \at 1. University College London, 2. The Chinese University of Hong Kong, Shenzhen 3. Department of Mathematical Sciences, University of Bath, 4.  Institut de Mathematiques de Toulouse, CNRS.
}

%
%
\maketitle
\abstract{We discuss parallel (additive) and sequential (multiplicative) variants of overlapping Schwarz methods for
    the Helmholtz equation in $\mathbb{R}^d$, with large real wavenumber and smooth variable
    wave speed. The radiation condition
is approximated by a Cartesian perfectly-matched layer (PML). The domain-decomposition
subdomains are overlapping hyperrectangles with Cartesian PMLs at their boundaries.  In a recent paper  ({\tt arXiv:2404.02156}),  the current  authors
proved (for both variants) that,    after a specified number of iterations -- depending on the behaviour of the geometric-optic rays -- the error is smooth and smaller than any negative power of the wavenumber $k$. For the parallel method, the specified number of iterations is less than the maximum number of subdomains, counted with their multiplicity, that a geometric-optic ray can intersect.  
  The theory, which is given at the continuous level and makes essential use of semi-classical analysis, assumes that the 
  overlaps of the subdomains and the widths of the PMLs are all 
   independent of
  the wavenumber. In this paper we extend the results of {\tt arXiv:2404.02156} by experimentally studying the behaviour of the methods   in the practically important case
  when both the overlap and the PML
  width decrease as the wavenumber increases. {We find that (at least for constant wavespeed), the methods remain robust to increasing $k$,  even for miminal  overlap,  when  the PML is one
    wavelength wide.}}

\section{The Helmholtz problem}

We consider the well-known Helmholtz equation:
\begin{equation}\label{eq:Helm}
- k^{-2} \Delta u - c^{-2} u = f \quad \text{in}\quad \mathbb{R}^d , 
\end{equation}
with the Sommerfeld radiation condition:
$r^{\frac{d-1}{2}} \left(\partial_r u - \mathrm{i} ku\right)\rightarrow 0, \
\text{as}\   r =|\mathbf{x}| \rightarrow \infty.$
Here  $k\geq 1$ is the wave number, $c \in C^\infty(\mathbb{R}^d)$ is the wavespeed and $f \in L^2(\mathbb{R}^d)$
is the source. 
While the paper \cite{GaGoGrLaSp:24} treats the case of general $d$, in the interests of brevity we restrict here to $d = 2$. The method for general $d$ is an obvious generalisation of the one presented here.
We  assume that  both $f$ and $1-c$ are  supported in a
box $\Omega_{\rm int}: = (0,l)\times(0,d).$

We now  restrict problem \eqref{eq:Helm} to the  domain
$\Omega := (-\kappa, l+\kappa)\times(-\kappa, d+\kappa)$, having added
a standard Cartesian PML of thickness   $\kappa$ to $\Omega_{\rm int}$.
To do this,  we choose a {\em scaling function}   $f_{s} \in C^{\infty}(\mathbb{R})$
(with $s$ denoting ``scaling''), satisfying 
$$
 f_{s}(x) = f'_{s}(x) = 0\text{ for }x\le 0\quad   \text{and} \quad f'_{s}(x)>0 \text{ for }x>0, \   ,
 $$
 together with  $f_{s}(x)'' = 0 \ \text{for} \ x \geq \kappa_{\rm lin}$
 for some $\kappa_{\rm lin} < \kappa$.
Using this,  we  define the horizontal and vertical scaling functions:   
$$
g_{1} (x_1) :=
\begin{cases}
f_{s}(x_1 - l ) &\text{if }x_1 \geq l  \\
0&\text{if }x_1 \in (0, l ) \\
{-} f_{s}( - x_1) &\text{if }x_1 \leq 0
\end{cases}, \hspace{0.3cm}
g_{2} (x_2) :=
\begin{cases}
f_{s}(x_2 - d ) &\text{if }x_2 \geq d  \\
0&\text{if }x_2 \in (0, d ) \\
{-} f_{s}( - x_2) &\text{if }x_2 \leq 0
\end{cases},
$$
and  the scaled operator 
\begin{align} \label{eq:scaledop}
\Delta_{s} : = \Big(
\gamma_{1}(x_1)^{-1} \partial_{x_1} \Big)^2 + \Big(
\gamma_{2}(x_2)^{-1} \partial_{x_2} \Big)^2, 
\end{align} 
where
$
\gamma_{m} := 1+ig_{m}'$, $m = 1,2$.
The PML approximation to \eqref{eq:Helm} then reads as: 
\begin{equation}\label{eq:PDE}
\begin{aligned}
\myP_s \, u  \ : = \ (-k^{-2} \Delta_{s}  - c^{-2})\,  u &= f \quad \text{on }\quad  \Omega,  \\
\text{subject to} \quad u&=0 \quad \text{on } \quad \partial \Omega. 
\end{aligned}
\end{equation}

\section{The overlapping Schwarz methods with local PMLs }
We consider both  parallel and sequential overlapping Schwarz methods for   \eqref{eq:PDE}, 
  using PML  as a subdomain boundary condition.  
Such a strategy was first proposed (without theory) in
\cite{To:98}, and has received much recent interest both in the overlapping and  non-overlapping cases --  for reviews see  \cite{GaGoGrLaSp:24} and   \cite{GaZh:22}.  
Recent work on the overlapping case includes \cite{LeJu:19} and  \cite{BoBoDoTo:22}. Our paper \cite{GaGoGrLaSp:24} provides the first wavenumber-explicit theory, 
both for the overlapping case,  and for any domain decomposition method for Helmholtz problems with a non-trivial scatterer.     

We cover  $\Omega_{\rm int}$
by  $N$ overlapping sub-rectangles: $\Omega_{{\rm int},j}: = (a^j,b^j)\times(c^j,d^j), \  j = 1,\ldots, N$, 
and we   let $\delta>0$ denote the minimum overlap parameter, as
defined for example in
\cite[Assumption 3.1]{ToWi:05}.
Each  $\Omega_{\rm int,j}$ is then extended  to a larger subdomain  $\Omega_j$ by adding a PML of width $\kappa$ along each edge,  analogous to the extension of $\Omega_{\rm int}$ to   $\Omega$.
(The PML width on interior edges could also be different from $\kappa$ - see \cite{GaGoGrLaSp:24}.)

An  example is the {\em checkerboard decomposition},  where the $\Omega_{{\rm int}, j} $ are constructed by  starting with  an $N_1 \times N_2$  tensor product rectangular (non-overlapping) decomposition of  $\Omega_{\rm int}$ and then extending   the internal
boundary of each sub-rectangle  outward,  subject to  the constraint that each extended sub-rectangle only overlaps
its nearest neighbours.  Then,  Cartesian PMLs are added to $\Omega_{\rm int}$ and $\Omega_{{\rm int},j}$ as above. A checkerboard with either $N_1 = 1$ or $N_2 = 1$ is called a
{\em strip decompositon}. In any case $N := N_1N_2$ is the number of subdomains.



We now define subproblems on $\Omega_{j}$, $j = 1, \ldots, N$.
The local scaling functions are  
$$
g_{1,j} (x_1) :=
\begin{cases}
f_{s}(x_1 - {b}_j ) &\text{if }x_1 \geq {b}_j  \\
0&\text{if }x_1 \in ({a}_j, {b}_j ) \\
{-} f_{s}( {a}_j- x_1) &\text{if }x_1 \leq {a}_j
\end{cases}, \hspace{0.3cm}
g_{2,j} (x_2) :=
\begin{cases}
f_{s}(x_2 - {d}_j ) &\text{if }x_2 \geq {d}_j \\
0&\text{if }x_2 \in ({c}_j, {d}_j ) \\
{-} f_{s}({c}_j - x_2) &\text{if }x_2 \leq {c}_j .
\end{cases}.
$$
Then the local scaled Laplace operators   (analogous to \eqref{eq:scaledop}) are   
\begin{align} \label{eq:localscaledop}
\Delta_{s,j} : = \Big(
\gamma_{1,j}(x_1)^{-1} \partial_{x_1} \Big)^2 + \Big(
\gamma_{2,j}(x_2)^{-1} \partial_{x_2} \Big)^2, 
\end{align} 
where
$
\gamma_{m,j} = 1+ig_{m,j}', m = 1,2$\ 
and the local scaled Helmholtz  operator 
is 
\begin{align}\label{eq:local_scaled}
P_{s,j} := - k^{-2} \Delta_{s,j}  - c^{-2} . 
\end{align}

To knit local solutions together, we introduce a 
partition of unity (POU)
on $\mathbb{R}^2$, denoted $\{ \chi_j\}_{j=1}^N$,
 such that   $\sum_j  \chi_j  \equiv 1  \ \text{on}\
  \mathbb{R}^2$,
and, in addition,   
\begin{align} \label{eq:POU}
 \supp \chi_j\cap \Omega  \subseteq \{ x \in \Omega_j: g_{\ell,j}(x_\ell)
  = g_\ell(x_\ell), \ \ell = 1,2 \}.
\end{align}
Condition \eqref{eq:POU} implies that  $\supp \chi_\ell \subseteq \Omega_{{\rm int},j}$ if $\Omega_{{\rm int},j} $ is an interior subdomain of  $\Omega_{\rm int}$
but  $\supp \chi_\ell$ is  extended into the PML of $\Omega_{\rm int}$ otherwise.  The algorithms are then as follows.

%

\medskip

\noindent 
{\bf  The additive (parallel) Schwarz method.}\   Given an initial guess  $u_+^0 \in H^1_0(\Omega)$,
we obtain the iterates $u_+^n: n = 1,2, \ldots$ by solving (variationally) the local problem:  
\begin{align}\label{eq:local_corrector}
P_{s,j}
\mathfrak c_j^n = 
  (f- P_{s}u_+^n )|_{\Omega_j} \quad \text{on} \quad \Omega_j  
  \end{align} for the corrector  $\mathfrak c_j^n \in H^1_0(\Omega_j)$,  and then setting 
\begin{align}\label{eq:adding_up}
u_+^{n+1} = u_+^n + \sum_{j=1}^N \chi_j \mathfrak c_j^n. \end{align} 

If we introduce  $u_j^{n+1} := u_+^n\vert_{\Omega_j} + \mathfrak c_j^n $,  then this can also be written:
\begin{align*}
P_{s,j}
u_j^{n+1} & = 
P_{s,j} (u_+^n|_{\Omega_j})
-(P_{s}u_+^n )|_{\Omega_j}
+ f|_{\Omega_j} \ \text{on} \ \Omega_j, \quad \text{with} \ u_j^{n+1} &= u_+^n  \quad \text{on} \ \partial \Omega_j, \\
\end{align*}
and then 
$
u_+^{n+1} := \sum_{j=1}^N \chi_j u_j^{n+1}$.
In this case, all the local contributions $u_j^{n+1}: j = 1,\ldots, N$ are computed independently in parallel, before $u_+^{n+1}$ is finally  assembled. The sequential version is a simple variant of this.     

\medskip

\noindent
{\bf The multiplicative (sequential) Schwarz method.} \   Given an initial guess $u_\times^0 \in H^1_0(\Omega)$, let $u^0_j:= u_\times^0|_{\Omega_j}.$  Then, for $n = 0, 1, \ldots$, do the following:  
 \begin{enumerate}
\item \label{it:for}(Forward sweeping) For $j=1,\ldots,N$, 
$$
u_{j,n}^{\rightarrow} := \sum_{\ell<j} \chi_\ell u_\ell^{2n+1} + \sum_{ \ell\geq j} \chi_\ell u_\ell^{2n},
$$
and then compute $u_j^{2n+1}\in H^1(\Omega_j)$ as  the solution to
\begin{align}\nonumber 
P_{s,j}
u_j^{2n+1} & = 
P_{s,j}(u_{j,n}^{\rightarrow}|_{\Omega_j})
-(P_{s}
u_{j,n}^{\rightarrow})|_{\Omega_j} + f|_{\Omega_j}  \quad \text{on} \quad \Omega_j, 
\\
\nonumber \text{subject to} \quad u_j^{2n+1} &= u_{j,n}^{\rightarrow}\quad \text{on} \quad \partial \Omega_j,
\end{align}

Then set
\beq\nonumber 
u_{\times}^{2n+1} := \sum_{\ell=1}^N \chi_\ell u_\ell^{2n+1}. 
\eeq
\item \label{it:back} (Backward sweeping)
  For $j=N,\ldots,1$, introduce 
$$
u_{j,n}^{\leftarrow} := \sum_{\ell\leq j} \chi_\ell u_\ell^{2n+1} + \sum_{j< \ell} \chi_\ell u_\ell^{2n + 2}\, \in H^1_0(\Omega).
$$
Then compute   $u_j^{2n+2} \in H^1(\Omega_j)$ as  the solution to
\begin{align}\nonumber 
P_{s,j}
u_j^{2n+2} & = 
P_{s,j}(u_{j,n}^{\leftarrow}|_{\Omega_j} ) 
-(P_{s}
             u_{j,n}^{\leftarrow} )|_{\Omega_j} + f|_{\Omega_j} \quad \text{on} \quad \Omega_j,
  \\
  \nonumber
\text{subject to} \quad u_j^{2n+2} &= u_{j,n}^{\leftarrow} \quad \text{on} \quad \partial \Omega_j, 
\end{align}
and then set  
\beq \nonumber 
u_{\times}^{2n+2} := u^{\leftarrow}_{N+1,n} = \sum_{\ell=1}^N \chi_\ell u_\ell^{2n+ 2}.
\eeq
\end{enumerate}

\begin{remark}
General sequential methods for any dimensional checkerboard decompositions can be found in \cite[Section 1.4.6]{GaGoGrLaSp:24}. These methods perform  multiple sweepings with different orders of the subdomains. Specifically, we construct exhaustive (see \cite[Definition 7.2]{GaGoGrLaSp:24}) sweeping methods such that, for each geometric-optic ray, there are at least two sweeps ordering the subdomains intersected along the ray in both forward and backward directions.
\end{remark}

\vspace{-1cm}
\section{Theoretical results}

Suppose the POU $\{\chi_j\}$ is $C^{\infty}$.  Let $\{u^n_*\}_{n=0}^\infty$ be any sequence of iterates for the additive ($*=+$) or multiplicative ($*=\times$) Schwarz algorithm above. Then  \cite[Theorems 1.1-1.4 and 1.6]{GaGoGrLaSp:24} give conditions that, for any  $M>0$ and integer $s \geq 1$,
guarantee the  existence of  $\cN \in \mathbb{N}$ and $C>0$
(both independent of $f$)  such that
\begin{align}\label{eq:conv}
  \Vert u - u^\cN_* \Vert_{H_k^s(\Omega)} \leq C k^{-M} \Vert u - u^0_* \Vert_{H_k^1(\Omega)}.\end{align}
(Here the norm is defined by
$
\N{v}^2_{H^s_k(\Omega)} := \sum_{|\alpha|\leq s} \N{ (k^{-1}\partial)^\alpha v}^2_{L^2(\Omega)}, \ \text{for} \ s \geq 1.
$)
 

In particular, \eqref{eq:conv} implies that the fixed-point iterations converge exponentially quickly in the number of iterations  for sufficiently-large $k$ and the rate of exponential convergence increases with $k$.

The case of no scatterer ($c \equiv 1$) is dealt with in \cite[Theorems 1.1-1.4]{GaGoGrLaSp:24}. In that case, for a strip DD  with $N$ subdomains, $\cN = N$ for the  parallel algorithm  while $\cN = 2$  for the sequential algorithm (one forward and backward sweep). For a 2D checkerboard with $N_1\geq2$ and $N_2 \geq 2$, we have $\cN = N_1 +N_2-1$ in  the parallel case  and $\cN = 4$ in  the sequential case (although the ordering of sweeps has to be carefully chosen).      The case of variable (but non-trapping) $c$ and general rectangular DD  is dealt with in  \cite[Theorems 1.6]{GaGoGrLaSp:24}, where  $\cN$ is defined as the maximum number of subdomains, counted with their multiplicity, that a geometric-optic ray can intersect.


These results are valid on fixed domains for sufficiently-large $k$, i.e., the PML widths and DD overlaps are arbitrary, but assumed independent of $k$.
Obtaining results that are also explicit in these geometric parameters of the decomposition will require more technical arguments than those used in \cite{GaGoGrLaSp:24}. In this paper we investigate this issue experimentally.  


\section{Discretisation and numerical results}

\noindent {\bf Discretisation.} \ Equation \eqref{eq:PDE},  involving the scaled operator $P_s$,  is recast in variation
form,  multiplying by  a test function $v \in H^1_0(\Omega)$ and integrating by parts to obtain
\begin{equation}\label{eq:global_vari}
a(u,v): = \int_{\Omega} k^{-2}\big((D\nabla u) \cdot \nabla \overline{v}  - (\beta\cdot \nabla u)\overline{v} \big)-  c^{-2}u \overline{v} = \langle f,  v\rangle_{L^2(\Omega)},    
\end{equation}
\beqs
\text{where} \quad D := \mathrm{diag}\left(  \gamma_{1}^{-2}(x_1), \gamma_{2}^{-2}(x_2)\right)
\
\tand\
\beta := \left(\gamma_{1}'(x_1)\gamma_{1}^{-3}(x_1)\, ,
\,  \gamma_{2}'(x_2) \gamma_{2}^{-3}(x_2) \right)^T . 
\eeqs
The variational form of equations involving the local operator $P_{s,j}$
(such as  \eqref{eq:local_corrector})
yield an analogous local sesquilinear form $a_j$ defined on $H^1_0(\Omega_j)^2$.


Let $\{\mathcal{T}_h\}$ be a shape-regular conforming sequence of meshes  for $\Omega$ which
resolve the boundaries of $\Omega$, $\Omega_{\rm int}$, $\Omega_{\rm int, j}$, and $\Omega_j$, for all $j$. Let $V_h$ be the space of continuous piecewise-polynomials
of degree $\leq p$ on $\mathcal{T}_h$ which vanish on $\partial \Omega$.
The finite element discretisation of  \eqref{eq:global_vari} leads to the linear algebraic system  
$\mathsf{A u =f}$, with $\mathsf{u}$ denoting the nodal vector of the finite element solution to be found. 

Then, with  $V_{h,j} = \{ v_h \vert_{\Omega_j}: v_h \in V_h\} \cap H^1_0(\Omega_j)$, 
the discrete version of the additive algorithm reads as follows.
Given an iterate $u^n_{+,h}\in V_h$, we compute a corrector $\mathfrak c_{h,j}^n \in V_{h,j}$ 
by solving the discrete version of  \eqref{eq:local_corrector}:
\beqs
a_j \big(\mathfrak c_{h,j}^n, v_{h,j} \big) = \big\langle f, \cR^T_{h,j} v_{h,j}\big\rangle - a\big(u^n_{h,j}, \cR^T_{h,j} v_{h,j} \big) \quad\tfa v_{h,j}\in V_{h,j},
\eeqs
where $\cR_{h,j}^\top$ denotes extension by zero from  $V_{h,j}$ to $V_h$. The  next iterate $u_{+,h}^{n+1}$ is obtained by the {discrete}  analogue of \eqref{eq:adding_up}. 

With $\mathsf{A}_j$ denoting the finite element stiffness matrix corresponding to $a_j$, and $\mathsf{u}_{+}^n$ denoting the nodal vector of the $n$th iterate, this can be easily seen to    correspond to the preconditioned Richardson iteration (see \cite[\S8]{GaGoGrLaSp:24}): 
\begin{align*}    \mathsf{u}^{n+1}_{+} := \mathsf{u}^n_{+} + \mathsf{B}^{-1} (f - \mathsf{A}
  \mathsf{u}^n_+)
  \quad \text{ with } \quad \mathsf{B}^{-1} :=  \sum_j \widetilde{\mathsf{R}}_{j}^\top
  \mathsf{A}_{j}^{-1}\mathsf{R}_{j},
\end{align*}
where $\mathsf{R}_j$ denotes the restriction of a nodal vector on $\Omega$ to its
nodes on $\Omega_j$, and $\widetilde{\mathsf{R}}_{j}^\top$ denotes the extension by zero
of a nodal vector on $\Omega_j$, after  multiplication by nodal values of $\chi_j$. Thus the additive algorithm is a familiar restricted additive Schwarz method where the
subdomain problems have PML boundary conditions.




\medskip

\noindent {\bf Numerical Experiment.} \ In this experiment, we {explore  how the overlap size $\delta$
and the definition of the PML affect}  the performance of the methods. For simplicity, we consider constant wavespeed $c=1$. 
{For discretization, we used finite elements of polynomial order $2$, with $h \sim k^{-1.25}$ (to control the pollution error).
  In the tests in \cite{GaGoGrLaSp:24}
  (chosen to illustrate the theory),  we had set    $\delta=\kappa=1/40$, fixed independently of $k$.
 In this case  the number of {freedoms} in the overlaps and the
PMLs  increases significantly as  $k$ increases, 
implying a heavy communication load at high-frequency.}
{Here  we investigate  smaller PML thickness $\kappa$ and smaller overlap   $\delta$, {thus increasing the practicality of the  method}.  The PML scaling function  is $f_{\rm PML}(x) = \frac{a}{3}x^3$.

{ We   tested two possible choices of  $\kappa$ and $a$. In the first setting  we chose  
  $\kappa = \lambda := \frac{2\pi}{k}$ (i.e., the PML is  one wavelength wide) and we compared   different choices of $a\in\{5\times10^3, \, 10^{4}, \, 30k, \, 0.3k^2\}$.
  In the second setting we chose  $\kappa \in \{ 3h, 5h\}$  and  $a = k^{2.5}$, (thus ensuring that the PML scaling functions have the same maximal modulus $\mathcal{O}(1)$ in the PML regions).
  The second setting 
  turned out  not to be    $k$-robust,  while the first was found to perform
  better 
 when  $a$ was increasing with
  respect to $k$. Due to the page limit, we only present the results for $\kappa = \lambda $ and $a = 30k$ and for overlap  $\delta \in \{ \frac{1}{80}, 2h, h\}$;
when $\delta = h$ the overlap has only two layers of elements.}

{Regarding the POU,  for a strip DD  we start from
functions $\tilde{\chi}_j$ supported on  $\Omega_{{\rm int}, j}$, which take value $1$ in the non-overlapping region and  decrease linearly  to $0$ towards the internal  boundary of $\Omega_{{\rm int}, j}$. Then the
POU is  obtained by normalizing, i.e.,  
$
\chi_j = \frac{\tilde{\chi}_j}{\sum_j\tilde{\chi}_j}.
$
For the checkerboard DD, we used the Cartesian product  $\tilde{\chi}_i(x)\tilde{\chi}_j(y)$ to generate local  functions on $\Omega_{{\rm int},j}$ and then normalized them to obtain $\chi_j$.
The POU used here is not as smooth as is required in the theoretical
analysis \cite{GaGoGrLaSp:24}. In fact, in  the experiments in \cite{GaGoGrLaSp:24}
we used  a smooth POU and each  $\chi_j$ was supported in a proper subset
of  $\Omega_{{\rm int}, j}$, bounded away from $\partial \Omega_{{\rm int}, j}$.   
}

Table \ref{tab:bigtable}  
lists the iteration counts {for the algorithms to attain} 1e-6 reduction of the relative residual. {Here  ``RAS-PML'' and ``RMS-PML''  denote respectively  the additive and 
  multiplicative algorithms,  each    tested on  strip and checkerboard DDs}.
For the RMS-PML {and}  checkerboard DD, we use the exhaustive sweeping method that contains 4 sweeps per iteration (see \cite[Figure 1.2]{GaGoGrLaSp:24}). 
  We  observe:
  1) The performance is hardly degraded  when the overlap is reduced from fixed  ($\delta = 1/80$) to minimal ($\delta = h$).
     When    $\delta =1/80$, the overlap contains    $4,7,10,13,16,19$ layers of
     elements for $k=100, 150, 200, 250, 300, 350$.
   The results show that this  communication cost 
   can be removed while hardly affecting the convergence.
2)
 With $\kappa = \lambda$, the iteration counts are not decreasing
    as fast with increasing $k$ as when we used  fixed  PML width. (See the
    experiments in \cite{GaGoGrLaSp:24} illustrating \eqref{eq:conv}.)
  However convergence rates remain bounded for the range of $k$ tested and
  in the  checkerboard case for RAS-PML,  the convergence rate decreases as the wave number increases; see Figure \ref{fig:convergence_rate}.


\begin{figure}[h!]
    \centering
    \begin{subfigure}[t]{0.5\textwidth}
        \centering
        \includegraphics[width=.85\textwidth]{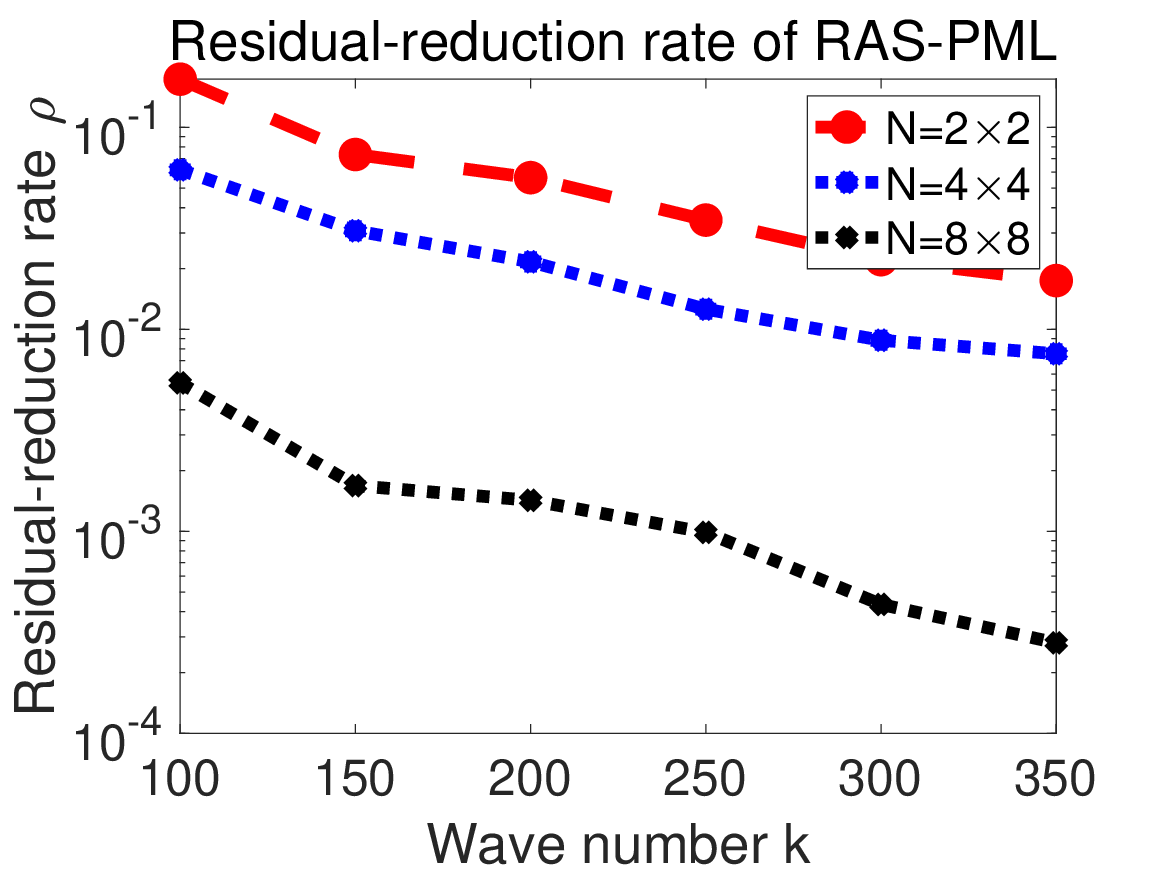}
        \caption{\footnotesize  RAS-PML: \tiny $\rho=\frac
        {\|\mathsf{A(u-u^{N_x+N_y-1})}\|_{\ell^2}}{\|\mathsf{A(u-u^0)}
        \|_{\ell^2}}$}
    \end{subfigure}%
         \hfill
    \begin{subfigure}[t]{0.5\textwidth}
        \centering
        \includegraphics[width=.85\textwidth]{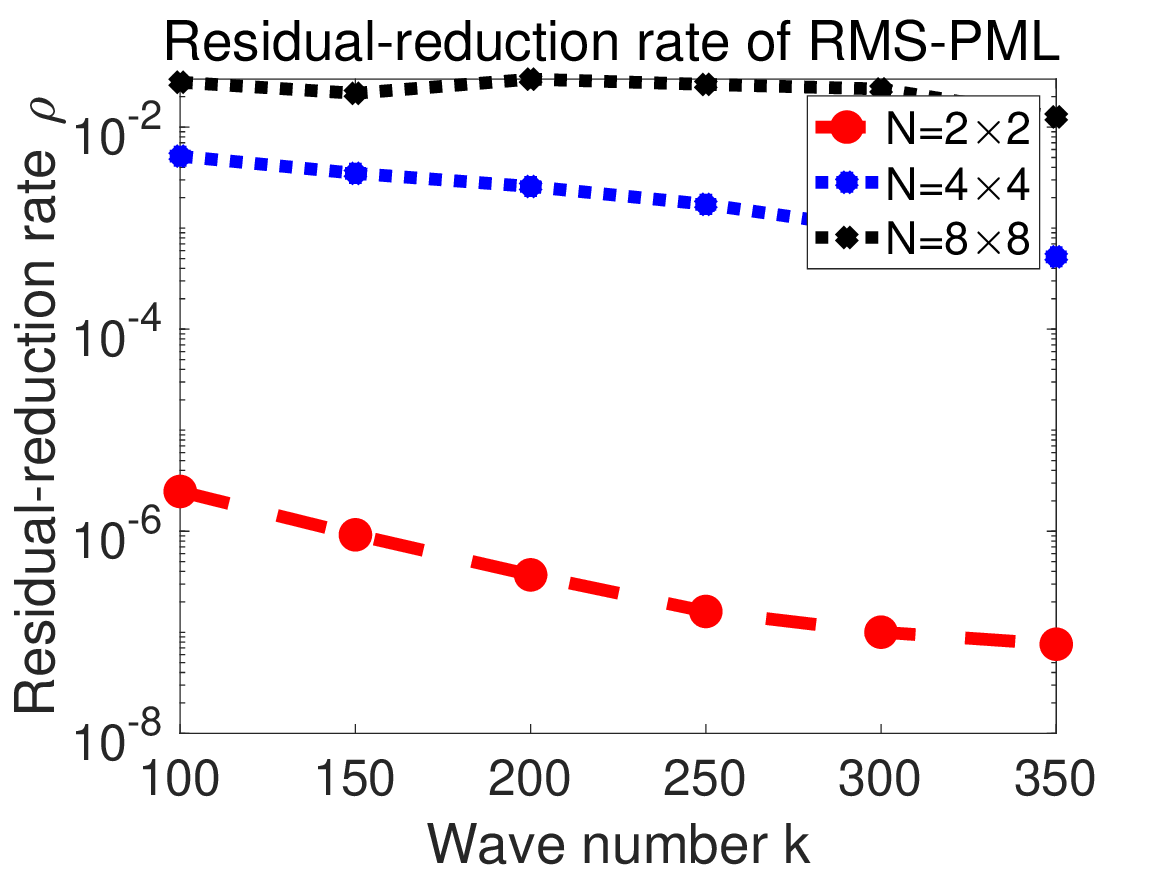}
        \caption{\footnotesize RMS-PML: \tiny$\rho=\frac
        {\|\mathsf{A(u-u^4)}\|_{\ell^2}}{\|\mathsf{A(u-u^0)}\|_{\ell^2}}$}
    \end{subfigure}%
   \caption{ RAS-PML/RMS-PML on Checkerboard: The rate of reduction of residual  against $k$. Left: Reduction rate per $N_x+N_y -1$ iterations for RAS-PML; Right: Reduction rate per exhaustive sweep for RMS-PML}\label{fig:convergence_rate}
\end{figure}

\begin{table}[h!]
  \centering
\begin{subtable}[h]{0.8\textwidth}     
  \centering
  \begin{tabular}{|c|ccc|ccc|ccc|}
\hline
N                         & \multicolumn{3}{c|}{2} & \multicolumn{3}{c|}{4} & \multicolumn{3}{c|}{8} \\ \hline
k\textbackslash{}$\delta$ & $\frac{1}{80}$    & 2h   & h   & $\frac{1}{80}$   & 2h   & h   & $\frac{1}{80}$   & 2h   & h  \\ \hline
100                        & 3(3)& 4(4)& 4(4)& 7(7)	& 7(7)	& 7(7)& 15(15)& 15(15)& 15(15)\\ 
150                        & 3(3)& 3(3)& 3(3)& 7(7)    	& 7(7)	& 7(7)& 15(15)& 15(15)& 15(15)\\
200                        & 3(3)& 4(3)& 4(3)& 7(7)	& 7(7)	& 7(7)& 12(12)& 15(15)& 15(15)\\
250                        & 4(4)& 4(4)& 4(4)& 6(6)	& 7(7)	& 7(7)& 11(11)& 15(15)& 15(15)\\ 
300                        & 4(4)& 5(5)& 5(5)& 7(7)	& 7(7)	& 7(7)& 15(15)& 15(15)& 15(15)\\
350                        & 4(4)& 5(5)&5(5)& 6(6)	& 7(7)	& 7(7)& 11(11)& 14(14)& 15(15)\\ \hline
 \hline
\end{tabular}
\caption{ RAS-PML on strip DD}
\end{subtable}
\label{tab:parallel1}
\begin{subtable}[h]{0.8\textwidth}     
\centering
\begin{tabular}{|c|ccc|ccc|ccc|}
\hline
N                         & \multicolumn{3}{c|}{$2\times 2$} & \multicolumn{3}{c|}{$4\times 4$} & \multicolumn{3}{c|}{$8\times 8$} \\ \hline
k\textbackslash{}$\delta$& $\frac{1}{80}$    & 2h   & h   & $\frac{1}{80}$    & 2h   & h   & $\frac{1}{80}$    & 2h   & h  \\ \hline
100                       & 8(6)& 9(6)& 9(7)&18(16)  	&19(18)& 20(19)& 33(27)& 40(32)& 47(36)\\
150                        & 8(6)& 8(6)& 8(7)& 16(15)	&18(18)&20(18)& 32(27)& 42(37)& 48(40)\\
200                        & 7(6)& 8(6)& 8(7)& 15(14)	& 18(17)& 19(18)& 30(28)& 48(37)& 56(40)\\
250                        & 7(6)& 8(6)& 8(6)& 14(13)	& 19(16)& 20(17)& 29(27)& 42(37)& 45(39)\\ 
300                        & 6(6)& 8(6)& 8(7)& 13(13)	& 18(16)& 18(17)& 28(27)& 44(37)& 47(39)\\
350                        & 6(6)& 7(6)& 7(7)& 13(13)	& 17(16)& 18(16)& 28(26)& 45(37)& 48(39)\\ \hline
\end{tabular}
\caption{ RAS-PML on checkerboard DD
}
\label{tab:parallel-checker1}

\end{subtable}

\begin{subtable}[h]{0.7\textwidth}     
\centering
    \begin{tabular}{lr}
\begin{tabular}{|c|ccc|ccc|ccc|}
\hline
N                         & \multicolumn{3}{c|}{2} & \multicolumn{3}{c|}{4} & \multicolumn{3}{c|}{8} \\ \hline
k\textbackslash{}$\delta$ & $\frac{1}{80}$    & 2h   & h   & $\frac{1}{80}$    & 2h   & h   & $\frac{1}{80}$    & 2h   & h  \\ \hline
100                        & 2& 2& 2& 2& 2& 2& 2& 2& 2\\
150                        & 2& 2& 2& 2& 2& 2& 2& 2& 2\\
200                        & 2& 2& 2& 2& 2& 2& 2& 2& 2\\
250                        & 2& 2& 2& 2& 2& 2& 2& 2& 2\\ 
300                        & 2& 3& 3& 2& 3& 3& 2& 3& 3\\
350                        & 2& 3& 3& 2& 3& 3& 2& 3& 3\\ \hline
\end{tabular}
\quad \quad & \quad \quad 
\begin{tabular}{|c|ccc|ccc|ccc|}
\hline
N                         & \multicolumn{3}{c|}{$2\times 2$} & \multicolumn{3}{c|}{$4\times 4$} & \multicolumn{3}{c|}{$8\times 8$} \\ \hline
k\textbackslash{}$\delta$ & $\frac{1}{80}$    & 2h   & h   & $\frac{1}{80}$    & 2h   & h   & $\frac{1}{80}$    & 2h   & h  \\ \hline
100                        & 2& 2& 2& 4& 4& 4& 5& 5& 6\\
150                        & 2& 2& 2& 3& 4& 4& 4& 5& 6\\
200                        & 2& 2& 2& 3& 4& 4& 4& 5& 6\\
250                        & 2& 2& 2& 3& 4& 4& 4& 5& 6\\ 
300                        & 2& 2& 2& 3& 3& 4& 4& 5& 5\\
350                        & 2& 2& 2& 3& 3& 3& 3& 5& 5\\ \hline
\end{tabular}
    \end{tabular}
         \caption{ RMS-PML: 
           \
 Left: strip DD, Right  checkerboard DD.}
\label{tab:sequential-checker1}

\end{subtable}

\caption{Iteration counts (relative residual tol = $10^{-6}$) with  $f_{\rm PML}(x) := \frac{a}{3}x^3$, \
  $a=30k$.  The numbers without bracket are for the fixed point iterations and the numbers in brackets are for the preconditioned GMRES iterations. Note that $\frac{1}{80h} = 4,7,10,13,16,19$ for $k=100,150,200,250,300,350$ \label{tab:bigtable}}
\end{table}




{\bf Acknowledgement} JG was supported by EPSRC grant EP/V001760/1. SG was supported by the National Natural Science Foundation of China (Grant number 12201535), the Guangdong Basic and
Applied Basic Research Foundation (Grant number 2023A1515011651) and  Shenzhen Stability Science Program 2022. DL was supported by
INSMI (CNRS) through a PEPS JCJC grant 2023. EAS was supported by EPSRC grant
EP/R005591/1. IGG was supported for several collaborative visits by CUHK Shenzhen.


\end{document}